\documentclass[aps,prl,preprint,groupedaddress]{revtex4-1}

\usepackage{graphicx}% Include figure files
\usepackage{dcolumn}% Align table columns on decimal point
\usepackage{bm}% bold math
\newcommand{\rmi}{\textrm{i}}
\newcommand{\rmd}{\textrm{d}}
\newcommand{\rme}{\textrm{e}}
\begin{document}

%\preprint{AIP/Chaos}

\title[An approach to normal forms of   Kuramoto model with distributed delays and the effect of minimal delay]{An approach to normal forms of   Kuramoto model with distributed delays and the effect of minimal delay}% Force line breaks with \\

\author{Ben Niu}

\email{niubenhit@163.com}%Lines break automatically or can be forced with \\
\affiliation{Department of  Mathematics, Harbin Institute of Technology(Weihai),\\  Weihai, 264200, P.R. China.}%

\author{Yuxiao Guo }
\affiliation{Department of  Mathematics, Harbin Institute of Technology(Weihai),\\  Weihai, 264200, P.R. China.}%
\author{Weihua Jiang }
%\email{niubenhit@163.com}%Lines break automatically or can be forced with \\

\affiliation{Department of  Mathematics, Harbin Institute of Technology,\\  Harbin 150001, P.R. China.}%

\date{\today}% It is always \today, today,

\begin{abstract}
 Heterogeneous delays with positive lower bound (gap) are taken into consideration in  Kuramoto oscillators. We first establish a perturbation technique, by which universal normal forms and detailed dynamical behavior of this model can be obtained easily.  Theoretically,   a hysteresis loop is found near the subcritically  bifurcated coherent state on the Ott-Antonsen's manifold.  For Gamma distributed delay with fixed variance and mean, we find large gap destroys the loop and significantly increases in the number of coexisted coherent attractors. This result is also explained in the viewpoint of excess kurtosis.
\end{abstract}

%\keywords{Kuramoto model; distributed delay with gap; bifurcation; normal form;synchronization}%Use showkeys class option if keyword
\pacs{}
\maketitle

\section{Introduction}

The Kuramoto phase oscillators were used to model diverse situations involving
large community of oscillators
\cite{article1001,article1002,article1003,article1004,article1005,article1006,article1007,article1008,article1009,article1110,article1111,article1112},
where the state of every oscillator is  determined by a
phase on the unit circle. This model captures essential features of
synchronization, observed in many physical models\cite{article1501,article1502,article1503,article1504,article1505} ranging from biology, neural science, lasers, engineering to superconducting
Josephson junctions.
Kuramoto\cite{article5} extended  Winfree's mean-field idea\cite{article10},
and confirmed that one population of weakly nearly-identical coupled
oscillators could be depicted as a universal model
$$
\dot\theta_i=\omega_i+\sum\limits_{j=1}^N\frac KN
\sin(\theta_j-\theta_i),~i=1,2,\ldots,N$$
Here the
 frequencies $\omega_i$ follow some distribution with
probability density function (PDF) $g(\omega)$.

In a network, signal's transmission and receive  both  lead to
 time delays.   Thus considering time lag is   necessary in many coupled systems\cite{ar9,ar92,ar921}. Due to  the spatio-distribution of oscillators, the transmitting delays among oscillators may be heterogenous. They may follow some probability distributions such as Gamma distribution, because the lag can be viewed as a period of awaiting. In a network with near-identical oscillators, the receiving delay, or responding delay, can be viewed as a constant. Hence, the total delay usually follows certain probability distribution with a gap, i.e., the positive lower bound. Other examples with a gap usually arise in the biomathematical problems. In population dynamics, the mature delay is an important parameter which is distributed  in an interval with positive lower bound\cite{gap1}. Thus in a prey-predator network, introducing heterogeneous delays with a gap should be greater realism. When dealing with different problems such as the growth of phytoplankton, the mature delay could also be replaced by time lag to digest nutrients\cite{gap2}. Time lags with a gap have also been used extensively when modelling  traffic flow dynamics\cite{gap3},  machine tool vibration problem\cite{gap4} and so on\cite{gap5}.

  Lee \textit{et al}\cite{wai}, using the Ott-Antonsen's
manifold reduction method\cite{ott2}, found that the
variation  of delay could greatly alter the dynamical behavior in a Kuramoto model with distributed delay without a gap.  This paper offered  a framework for  studying the  delay heterogeneity, where the results are illustrated with
respect to Gamma-distributed time lags in $[0,+\infty)$. After some simulations,
both supercritical and subcritical Hopf bifurcations on the
mean-field are observed.

\begin{figure}[htbp]
 (a)~~~~~~~~~~~~~~~~~~~~~~~~~~~~~~~~~~~~~~~~~~~~~~~~~~~~~~~~~~~(b)\\
  \includegraphics[width=0.45\textwidth]{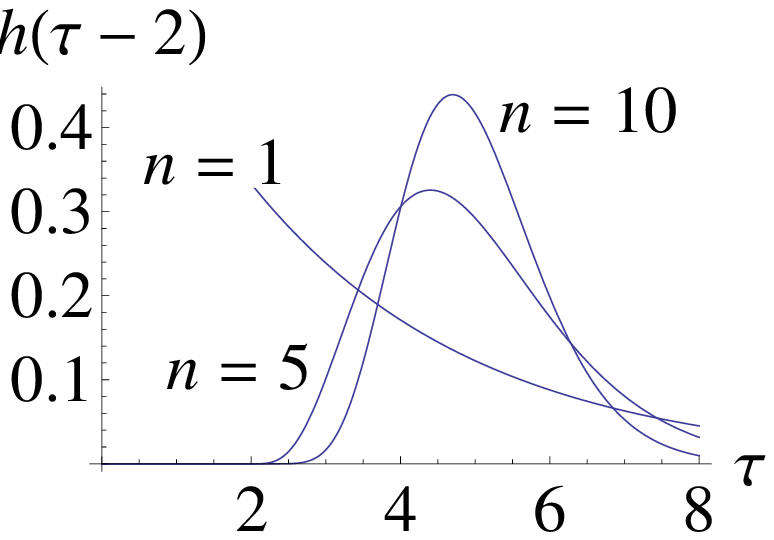}~
 \includegraphics[width=0.45\textwidth]{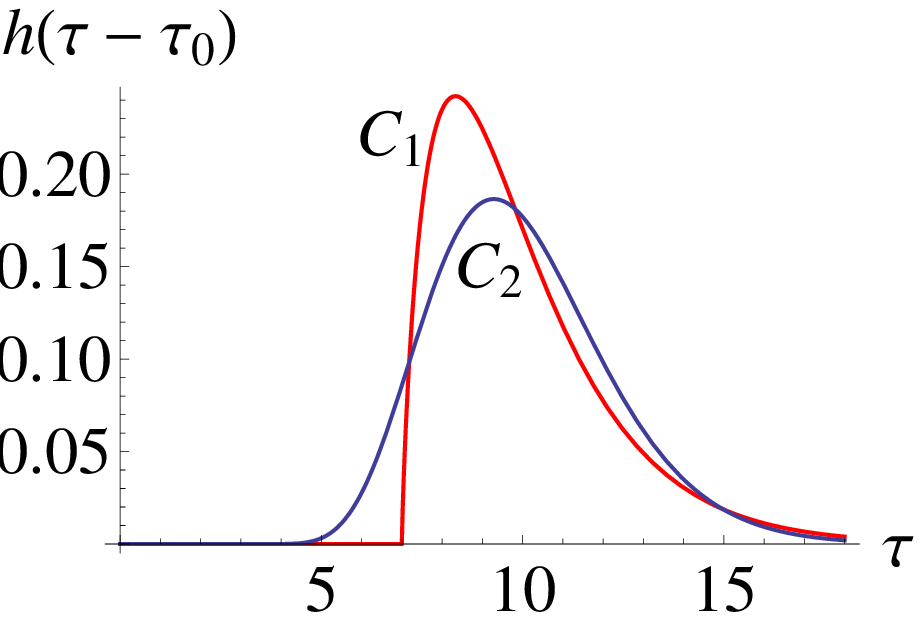}
  \caption{(a) $h(\tau-\tau_0)$,  PDFs of the sum of Gamma distribution
$\tilde\tau\sim\Gamma(n,\frac T n)$ with $T=3$, $n=1,5,10$ and $\tau_0=2$. (b) Two PDFs of Gamma distribution with a gap, $C_1$: (red) $T=3,\tau_0=7$ with larger excess kurtosis; $C_2$: (blue) $T=7,\tau_0=3$ with smaller excess kurtosis. $n$ is suitably   chosen such that variance $T^2/n=5$.
}\label{gammafig}
\end{figure}

 On one hand, in mathematical consideration, the mechanism causing the above phenomena is not quite clear yet, which depends on further bifurcation analysis. Using bifurcation technique, the transition among different schemes can be detected clearly.   On the other, when    time lag distributes in the interval $[\tau_0,+\infty)$ with a minimal responding time $\tau_0>0$, this is the case rarely investigated before. The total delay may be the sum of a Gamma distributed delay $\tilde\tau\sim\Gamma(n,\frac T n)$  and a constant  $\tau_0$ (See FIG. \ref{gammafig}(a)), for example.  In this case, an interesting fact is that studying only the mean and variance of delays sometimes does not make any senses.    For instance,  by varying $n$, one can still fix the  expectation of total delay $\langle\tau\rangle=\langle\tilde\tau+\tau_0\rangle=T+\tau_0$, and its variance $\textrm{Var}(\tau)=\frac {T^2} n$ despite the ratio $\frac{\tau_0}{T+\tau_0}$ varies (See FIG. \ref{gammafig}(b)). Thus the gap $\tau_0$ may have certain effect on the system dynamics without changing $\langle\tau\rangle$ and $\textrm{Var}(\tau)$. In this case higher order  moments of the data should be considered such as skewness or excess kurtosis, which involve the third-order or fourth-order central moments and are usually used to measure  the ``asymmetry'' or the  ``peakedness'' of  probability distribution, respectively. So far as we know, these two points of view are new and have not been well studied.

Motivated by the above two considerations and these pioneer works, we are about to consider a more realistic case: heterogeneous  delays with a positive minimal delay in a system, or  a gap\cite{yuan}, we say.
Now the Kuramoto model reads
\begin{equation}\label{model}
\dot
\theta_i=\omega_i+\frac{k}{N}\sum_{j=1}^{N}\sin[\theta_j(t-\tau_i)-\theta_i(t)]\end{equation}
$i=1,2,\ldots,N$,
where $\theta_i(t)\in[0,2\pi)$ is the phase of the $i$th oscillator
and $\omega_i$ is its natural frequency coming from an ensemble with PDF $g(\omega), \omega\in(-\infty,+\infty)$. $k$ is
the constant coupling strength. The
delays $\tau_i$, coming from the ensemble $\tau$, are   statistically
independent with $\omega_i$.  In  this paper, we choose random variable $\tau=\tilde \tau+\tau_0$ with $\tilde\tau$ having PDF $h(\tau)$ and $\tau_0$ a constant. Thus $\tau_i$ can be viewed as following a new distribution with PDF $h(\tau-\tau_0), \tau\in[\tau_0,+\infty)$.

In this paper, we will give some bifurcation results about system (\ref{model}).  After reducing it onto the Ott-Antonsen's manifold, a delay differential equation is obtained, in which   a Hopf bifurcation at the trivial solution means the coherent state is bifurcating  from the incoherent state.  In our earlier work\cite{niuphyd}, the center manifold reduction method\cite{Hassardfaria,Hassardfaria} is employed to investigate this bifurcation. However, the calculations depend on a rather complicated decomposition of a Banach space and many Riemann-Stieltjes integrals. Here, we make this approach much easier, and use the method of multiple scales\cite{nams,nams1,nams2} to give a relatively simple calculation process of the normal forms, by which the direction and the stability of the bifurcated coherent state are determined. For Gamma distributed delay, we   calculate all bifurcation points in certain parameter spaces and discuss the effect of the gap. Finally, it is found that, when fixing $\langle\tau\rangle$ and Var$(\tau)$, larger gap (or larger excess kurtosis, equivalently) not only leads to a supercritical bifurcation hence avoids the existence of hysteresis loop, but also significantly increases in the number of coexisted coherent states.

\section{Reduction}

As $N\rightarrow\infty$, the continuity equation of    (\ref{model}) is
\begin{equation}\label{fp}
\frac\partial{\partial t}f+\frac\partial{\partial \theta}(\zeta
f)=0\end{equation}with  a drift term $\zeta(\theta,t)=\omega+\frac k{2\rmi}(\rme^{-\rmi\theta}r-\rme^{\rmi\theta}r^\ast)$. The complex-valued ``order-parameter'' $r(t)$ is defined by
\begin{equation}\label{r}
r(t)=\langle\xi(t-\tau)\rangle=\int_{\tau_0}^{+\infty}\xi(t-\tau)h(\tau-\tau_0)\rmd\tau
\end{equation}
with
\begin{equation}\label{xi}
\xi(t)=\int_{-\infty}^{+\infty}\int_0^{2\pi}f(\omega,\theta,t)\rme^{\rmi\theta}\rmd\theta \rmd\omega
\end{equation}
The distribution density $f(\omega,\theta,t)$ characterizes the state of the oscillators' system at time $t$ in frequency $\omega$ and phase $\theta$.

 Now we are about to restate some results about the Ott-Antonsen's reduction of a system with distributed delay first derived by Lee \textit{et al}\cite{wai}.
Rewriting system (\ref{fp}) as
\begin{equation}\label{fokkerplanck}
\frac\partial{\partial t}f+\frac\partial{\partial
\theta}\left\{\left[\omega+\frac
k{2\rmi}(\rme^{-\rmi\theta}r-\rme^{\rmi\theta}r^\ast)\right]
f\right\}=0\end{equation} and restricting this partial differential
equation on the Ott-Antonsen manifold
$$
\left\{f:f=\frac{g(\omega)}{2\pi}\left\{1+\left[\sum_{m=1}^\infty\alpha^m(\omega,t)\rme^{\rmi m\theta}+c.c.\right]\right\}\right\}$$
with c.c. the complex conjugate of the formal terms,  we substitute the Fourier series  of $f$ into  (\ref{fokkerplanck}). After comparing the coefficient of the  same harmonic terms, a reduced equation is obtained
\begin{equation}\label{alpha}
\dot \alpha(\omega,t)=-\rmi \omega\alpha(\omega,t)+\frac k 2 r^\ast-\frac k 2 r \alpha^2(\omega,t)
\end{equation}
Obviously, from (\ref{xi}) we have $\xi(t)=\int_{-\infty}^{+\infty}g(\omega)\alpha^\ast(\omega,t)\rmd\omega$, then Eq.(\ref{r}) yields
\begin{equation}\label{rt}
r(t)=\int_{\tau_0}^{+\infty}\int_{-\infty}^{+\infty}g(\omega)\alpha^\ast(\omega,t-\tau)\rmd\omega h(\tau-\tau_0)\rmd\tau\end{equation}

For the sake of theoretical analysis, the distribution density
$g(\omega)$ is usually chosen as  Lorentzian distribution, that is
\begin{equation}\label{loren}
g(\omega)=\frac{\Delta}{\pi[(\omega-\omega_0)^2+\Delta^2]},~-\infty<\omega<+\infty\end{equation}
Note that the Lorentzian distribution is unimodal and can be viewed as an approximation of normal distribution.

Following Ott and Antonsen's method \cite{ott2},  substituting (\ref{loren}) into (\ref{rt}), and using residue theorem, we have
\begin{equation}\label{sys1}
r(t)=\int_{\tau_0}^{+\infty}\alpha^\ast(\omega_0-\rmi\Delta,t-\tau)h(\tau-\tau_0)\rmd\tau\end{equation}
Putting $\omega=\omega_0-\rmi\Delta$ in Eq.(\ref{alpha}) and
noticing Eq.(\ref{sys1}) yield
\begin{equation}\label{sys2}\begin{array}{l}
\dot \alpha(t)=-(\rmi\omega_0+\Delta)\alpha(t)+\frac k 2 \int_{\tau_0}^{+\infty}\alpha(t-\tau)h(\tau-\tau_0)\rmd\tau
-\frac k 2 \int_{\tau_0}^{+\infty}\alpha^\ast(t-\tau)h(\tau-\tau_0)\rmd\tau \alpha^2(t)\end{array}
\end{equation}
which is a delay differential equation\cite{Halefde}, whose trivial equilibrium stands for the incoherent state. To investigate its stability, we substitute
$\alpha=\alpha_0\rme^{\lambda t}$ into   the linear part of (\ref{sys2}) with $\alpha_0\neq0$, and obtain
\begin{equation}\label{CE}
\lambda=-(\rmi\omega_0+\Delta)+\frac k 2
\int_{\tau_0}^{+\infty}\rme^{-\lambda\tau}h(\tau-\tau_0)\rmd\tau\end{equation}

Obviously, if $\Delta>0$, $\lambda=-(\rmi\omega_0+\Delta)$ at $k=0$, thus the incoherence is stable for sufficiently small $k$. After the occurrence of Hopf bifurcation, roots with positive real part may appear in (\ref{CE}), which means the incoherence looses stability.   At the  limit of identical oscillators $\Delta\rightarrow 0$, we know the incoherence is neutrally stable at $k=0$. Once $k$ increases, the incoherent state becomes stable (or unstable) if Re~$\frac{\rmd\lambda}{\rmd k}|_{k=0}=\frac 12 \int_{\tau_0}^{+\infty}\cos\omega_0\tau h(\tau-\tau_0)\rmd\tau<0$ (or $>0$).

Usually, if (\ref{CE}) has a root $\lambda=\rmi\beta$,  a Hopf bifurcating solution $\alpha(t)=\rme^{\rmi\beta t}$ appears, which stands for a coherent state in system (\ref{model}). In the rest part of this paper, the stability of this coherent state and the exact location where it appears will be investigated.

\section{Bifurcation analysis}

In this section, we assume that $\Delta>0$, and  a Hopf bifurcation occurs in
Eq.(\ref{sys2}).
If Hopf bifurcation occurs at $k=\bar{k}$, two necessary conditions are required:
 [i] Eq.(\ref{CE}) has a simple root $\lambda=\rmi\beta$ with $\beta\neq 0$
  when $k=\bar{k} $;
 [ii] The so-called  transversality condition holds in the sense that $\textrm{Re}\lambda'(k)\neq0$ at $k=\bar{k} $.

 To obtain more properties near $\bar{k} $, one usually should do bifurcation analysis near this critical point, including the normal form deriving  and unfolding analysis\cite{Wiggins2}. In a previous work\cite{niuphyd}, we have done this in case of $\tau_0=0$. The approach therein heavily depends on a mathematical fundation such as formal adjoint theory of functional differential equation and the decomposition of Banach space. Moreover, some calculations such as Riemann-Stieltjes integrals are tedious, and we can expect that even more calculations are needed when $\tau_0>0$. In this section, we are about to extend the traditional method of multiple scales\cite{nams}, a classical perturbation method, to the bifurcation analysis of the complex-valued equation (\ref{sys2}).

\subsection{Normal forms}

Denoting by $k=\bar{k} +\epsilon\nu$ with $\epsilon>0$ and  $\nu$ a detuning parameter which describes the nearness of $k$ to the critical value  $\bar k$, Eq.(\ref{sys2})
can be rewritten into
\begin{equation}\label{sys}\begin{array}{ll}
\dot \alpha(t)=&-(\rmi\omega_0+\Delta)\alpha(t)+\frac
{\bar{k} +\epsilon\nu} 2
\int_{\tau_0}^{+\infty}\alpha(t-\tau)h(\tau-{\tau_0})\rmd\tau
\\&-\frac
{\bar{k} +\epsilon\nu} 2
\int_{\tau_0}^{+\infty}\alpha^\ast(t-\tau)h(\tau-{\tau_0})\rmd\tau
\alpha^2(t)\end{array}
\end{equation}

\begin{figure}[htbp]
   (a)~~~~~~~~~~~~~~~~~~~~~~~~~~~~~~~~~~~~~~~~~~~~~~~~~~~~~~~~~~~(b)\\
  \includegraphics[width=0.48\textwidth]{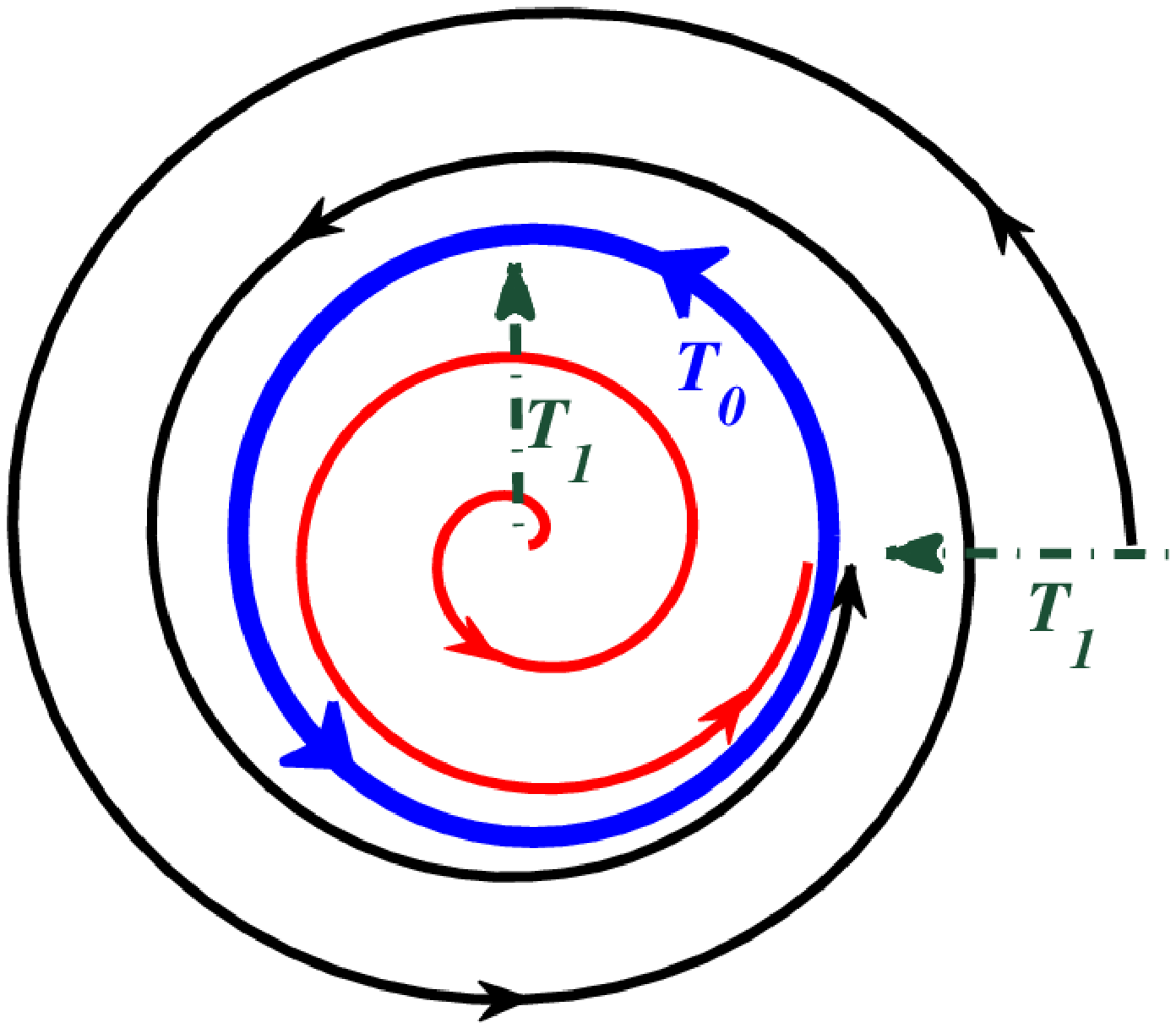}
   \includegraphics[width=0.48\textwidth]{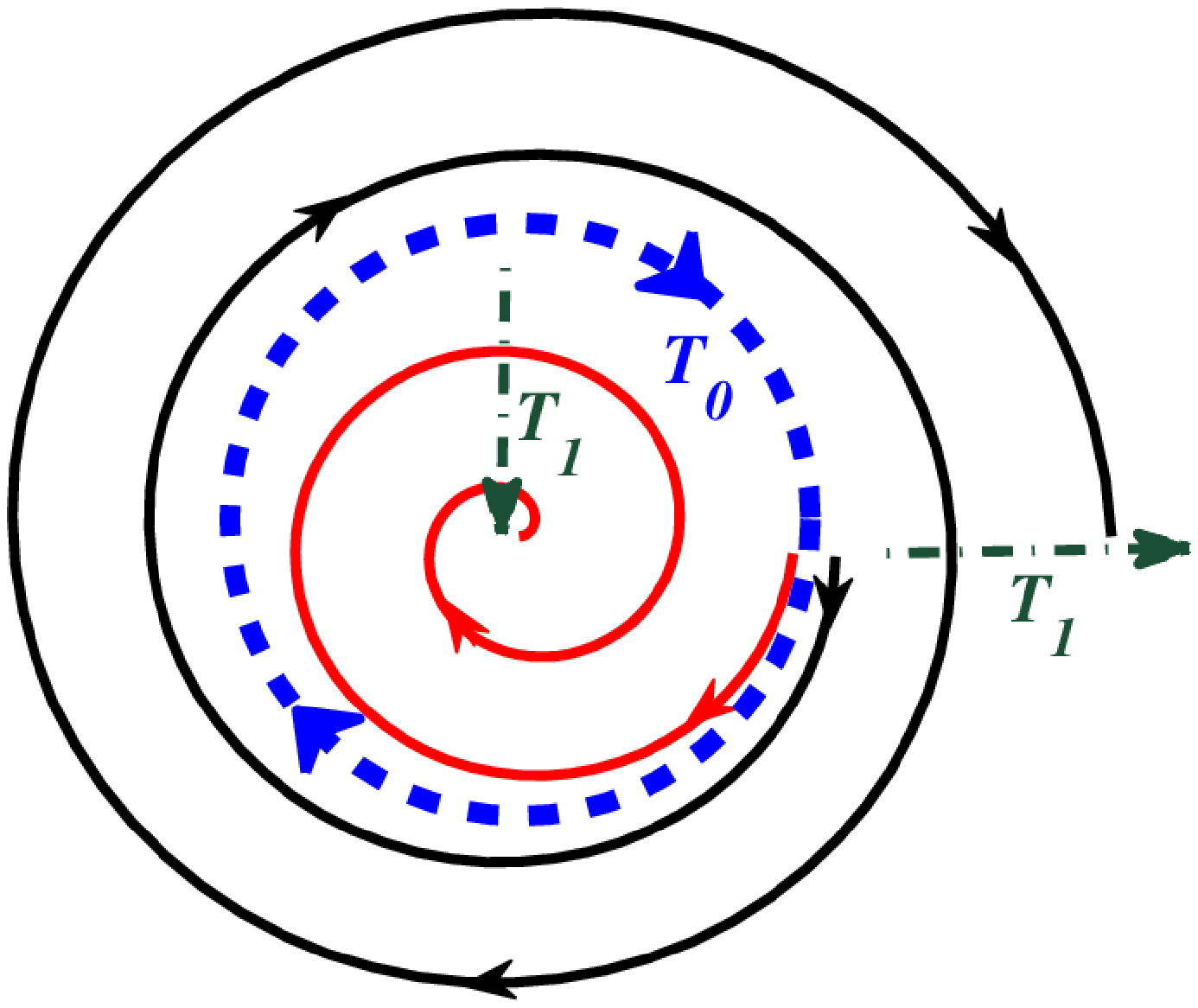}
  \caption{(a) Solutions  trend towards (in scale $T_1$, red and black curves) a stable cycle (oscillating in scale $T_0$, blue).
  (b) Solutions  depart from (in scale $T_1$, red and black curves) an unstable cycle (oscillating in scale $T_0$, blue).
}\label{appro}
\end{figure}

  For the
absence of second order term in (\ref{sys}), the solution to
Eq.(\ref{sys}) can be expressed by\cite{epsilon22}
$$
\alpha(t;\epsilon):=\epsilon^{1/2}\alpha_1(T_0,T_1)+\epsilon^{3/2}\alpha_2(T_0,T_1)+\cdots$$
where $T_0=t$, $T_1=\epsilon t$. The two time scales can be interpreted as follows. When $\nu$ is sufficiently small, i.e., near the critical point,  solutions to system (\ref{sys}) near the trivial equilibrium oscillate in a fast time scale $T_0$\cite{explain}, whereas they trend towards (\textit{OR} depart from) a small stable (\textit{OR} unstable) periodic oscillation in a slow time scale $T_1$. This is shown in FIG. \ref{appro}.

 The derivative with respect to $t$
is transformed into
$$
\frac{\rmd}{\rmd t}=\frac{\partial}{\partial
T_0}+\epsilon\frac{\partial}{\partial T_1}:=D_0+\epsilon
D_1$$
Taylor expanding the term $\alpha(t-\tau):=\alpha(t-\tau;\epsilon)$ gives
$$\begin{array}{l}
\alpha(t-\tau;\epsilon)=\epsilon^{1/2}\alpha_1(T_0-\tau,T_1)+\epsilon^{3/2}\alpha_2(T_0-\tau,T_1)-\epsilon^{3/2}\tau D_1\alpha_1(T_0-\tau,T_1)+\cdots\end{array}$$

Substituting $\alpha(t)$ and $\alpha(t-\tau)$ into (\ref{sys}), and balancing the same order terms of $\epsilon$ in both sides we have
\begin{equation}\label{ep12}\begin{array}{l}
D_0 \alpha_1(T_0,T_1)=-(\rmi\omega_0+\Delta)\alpha_1(T_0,T_1)+\frac
{\bar{k} } 2
\int_{\tau_0}^{+\infty}\alpha_1(T_0-\tau,T_1)h(\tau-{\tau_0})\rmd\tau\end{array}
\end{equation}
and
\begin{equation}\label{ep13}\begin{array}{ll}
D_0 \alpha_2(T_0,T_1)=&-(\rmi\omega_0+\Delta)\alpha_2(T_0,T_1)+\frac
{\bar{k} } 2
\int_{\tau_0}^{+\infty}\alpha_2(T_0-\tau,T_1)h(\tau-{\tau_0})\rmd\tau\\&
+\nu \frac
{1} 2
\int_{\tau_0}^{+\infty}\alpha_1(T_0-\tau,T_1)h(\tau-{\tau_0})\rmd\tau-D_1 \alpha_1(T_0,T_1)\\&+ \frac
{\bar{k} } 2
\int_{\tau_0}^{+\infty}-\tau D_1\alpha_1(T_0-\tau,T_1)h(\tau-{\tau_0})\rmd\tau\\&-\frac
{\bar{k} } 2
\int_{\tau_0}^{+\infty} \alpha_1^\ast(T_0-\tau,T_1)h(\tau-{\tau_0})\rmd\tau\alpha_1^2(T_0,T_1)\end{array}
\end{equation}
 Eq.(\ref{ep12}) is a linear equation and has a solution
$$\begin{array}{l}
\alpha_1(T_0,T_1)=A(T_1)\rme^{\rmi \beta T_0}+\sum\limits_{m=1}^{+\infty}[A_m(T_1)\rme^{(\sigma_m+\rmi
\beta_m)   T_0}]\end{array}
$$
If all  roots of Eq.(\ref{CE}) $\sigma_m+\rmi
\beta_m,~m=1,2,\ldots$, except $\rmi\beta$, have negative real part, $\sigma_m<0$ for any $m$, then $\alpha_1(T_0,T_1)\rightarrow A(T_1)\rme^{\rmi \beta T_0}$, as $t\rightarrow+\infty$.

 Now we can find that, the bifurcated solution oscillates in time scale $T_0$, and all solutions nearby trend towards it in time scale $T_1$, if the periodic oscillation is locally attractive. Thus the dynamical behavior of $A(T_1)$ determines the property of the bifurcation such as the stability: if a differential equation about $|A(T_1)|$ has a stable nontrivial equilibrium, system (\ref{sys}) has a stable periodic solution originating from Hopf bifurcation. This is the main idea of normal form method of Hopf bifurcation [recalling FIG. \ref{appro}],   which will be further calculated in the following.

Substituting $\alpha_1(T_0,T_1)$ into (\ref{ep13}) yields
$$\begin{array}{ll}
D_0 \alpha_2(T_0,T_1)&=-(\rmi\omega_0+\Delta)\alpha_2(T_0,T_1)+\frac
{\bar{k} } 2
\int_{\tau_0}^{+\infty}\alpha_2(T_0-\tau,T_1)h(\tau-{\tau_0})\rmd\tau\\&+\nu \frac
{1 } 2
\int_{\tau_0}^{+\infty} A(T_1)\rme^{\rmi \beta (T_0-\tau)}h(\tau-{\tau_0})\rmd\tau
- A'(T_1)\rme^{\rmi \beta T_0}\\&+ \frac
{\bar{k} } 2
\int_{\tau_0}^{+\infty}-\tau  A'(T_1)\rme^{\rmi \beta (T_0-\tau)}h(\tau-{\tau_0})\rmd\tau\\&-\frac
{\bar{k} } 2
\int_{\tau_0}^{+\infty}  A^\ast(T_1)\rme^{-\rmi \beta (T_0-\tau)}h(\tau-{\tau_0})\rmd\tau A^2(T_1)\rme^{2\rmi \beta T_0}\end{array}
$$

The last four terms in the above equation lead to secular terms, because they make $\alpha_2(T_0,T_1)$ depend on $T_0\rme^{\rmi \beta T_0}$, which is a contradiction.  Eliminating them, we have the normal form given by
$$\begin{array}{ll}
A'(T_1)&=\nu \frac{
\int_{\tau_0}^{+\infty} \rme^{-\rmi \beta \tau}h(\tau-{\tau_0})\rmd\tau}{2+
{\bar{k} }
\int_{\tau_0}^{+\infty}\tau   \rme^{-\rmi \beta  \tau}h(\tau-{\tau_0})\rmd\tau}A(T_1)+\frac{
{\bar{k} }
\int_{\tau_0}^{+\infty}  \rme^{\rmi \beta \tau}h(\tau-{\tau_0})\rmd\tau}{2+
{\bar{k} }
\int_{\tau_0}^{+\infty}\tau   \rme^{-\rmi \beta  \tau}h(\tau-{\tau_0})\rmd\tau}A^2(T_1)A^\ast(T_1) \\&:=a\nu A(T_1)+bA^2(T_1)A^\ast(T_1)\end{array}
$$

 The amplitude equation is
$$|A(T_1)|'=\nu\textrm{Re} a  |A(T_1)|+\textrm{Re} b |A(T_1)|^3$$

In fact, regarding $\lambda$ a function of $k$ defined implicitly by (\ref{CE}), we know $a=\lambda'(\bar k)$. According to the fundamental theory about Poincar\'e--Birkhoff normal form of ODE\cite{Wiggins2}, the two real parts of $a$ and $b$ determine the direction and stability of the bifurcation.

Precisely, letting
$k_0=\inf\{\bar{k} |\bar{k} $ is a Hopf bifurcation value$\}$,
we  always have Re$a>0$. If $\nu>0$ and Re$b<0$, then $|A(T_1)|\rightarrow \sqrt{\frac{\nu\textrm{Re}a}{-\textrm{Re}b}}$ as $T_1\rightarrow\infty$. Thus a branch of stable\cite{stable} bifurcation solutions appears at $k>k_0$ for Re$b<0$. Similarly, a branch of  unstable solutions appears at $k<k_0$ for Re$b>0$. In the latter case, the branch of unstable solutions must go to $k>k_0$, because when $k=0$ system (\ref{sys2}) has a globally stable attractor $\alpha=0$.

\begin{figure}[htbp]
  \centering\includegraphics[width=0.9\textwidth]{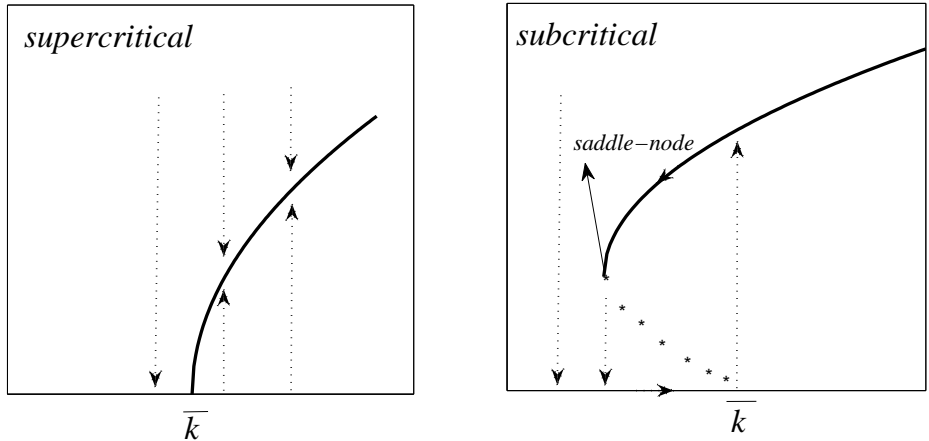}
  \caption{Supercritical bifurcation and subcritical bifurcation(with a hysteresis loop between $\bar k$ and the saddle-node point) near a critical value $\bar k$.}\label{exa}
\end{figure}

Two kinds of bifurcations are the so-called supercritical and subcritical bifurcations, respectively as shown
in FIG. \ref{exa}. If at some $\bar k$, Re $a<0$, there are two other kinds of bifurcations which can be illustrated by reversing the stability of coherent states.

Consider the case  time delay comes from an ensemble $\tau$ which is the sum of
a Gamma distributed variable $\tilde \tau$ and a constant $\tau_0$. $\tilde\tau\sim\Gamma(n,\frac T n)$ has PDF
$$
h(\tau)=\frac 1
{\Gamma(n)(T/n)^n}\tau^{n-1}\rme^{-\frac\tau{T/n}},\tau>0$$
where   $T$ is the   mean value, $\frac {T^2} n$ the variance, $\frac{2}{\sqrt n}$ the skewness and $\frac 6 n$ the excess kurtosis. The PDF of $\tau$, $h(\tau-\tau_0)$ has been shown in FIG. \ref{gammafig}(a). The laplace transformation of $h(\tau)$ is $(1+\frac T
n\lambda)^{-n}$.

The characteristic equation in this case is given by
$$\begin{array}{l}
\lambda=-(\rmi\omega_0+\Delta)+\frac k
{2}\int_{\tau_0}^{+\infty}\rme^{-\lambda
\tau}h(\tau-{\tau_0})\rmd\tau=-(\rmi\omega_0+\Delta)+\frac k {2}(1+\frac T
n\lambda)^{-n}\rme^{-\lambda\tau_0}\end{array}
$$
Motivated by the above analysis, we are seeking for a root $\lambda=\rmi\beta$ and
have
\begin{equation}\label{cegammas}
\rmi\beta=-(\rmi\omega_0+\Delta)+\frac k {2}(1+\frac T
n\rmi\beta)^{-n}\rme^{-\rmi\beta\tau_0}
\end{equation}
Denote by $1+\frac T n\rmi\beta=\rho\rme^{\rmi\vartheta}$, then
\begin{equation}\label{cossingamma}\begin{array}{rl}
\Delta=&\frac k 2 \rho^{-n}\cos (n\vartheta+\beta\tau_0)\\
\beta+\omega_0=&-\frac k 2 \rho^{-n}\sin
(n\vartheta+\beta\tau_0)\end{array}\end{equation}
Obviously, $\vartheta\in[-\pi/2,\pi/2]$, thus $\vartheta=\arctan \frac{T
   \beta }{n}$. Then we have$$-\frac{\beta
+\omega _0}{\Delta }=\tan \left[n \arctan\left(\frac{T
   \beta }{n}\right)+\beta\tau_0\right]$$
This equation can be solved by a sequence of $\beta$'s with $|\beta|\rightarrow+\infty$, for all parameters fixed.
As $|\beta|\rightarrow+\infty$, we have $\rho\rightarrow\infty$, thus $k\rightarrow\infty$ holds from the first equation of (\ref{cossingamma}).
Hence we only need to calculate all roots for $\beta\in[-J,J]$ with $J$ relatively large, to obtain the first several bifurcation values of $\bar k$.

If a  bifurcation occurs at $\bar k$, then employing  the normal form theory established above, we have
$$a=\frac{(1+\frac {\rmi\beta T}
n)^{-n}\rme^{-\rmi\beta\tau_0}}{2+\bar k \left(\frac{T \rme^{-\rmi \beta  \tau_0 }}{ 1+\frac{\rmi
   \beta  T}{n} }+\tau_0  \rme^{-\rmi \beta  \tau_0 }\right) \left(1+\frac{\rmi
   \beta  T}{n}\right)^{-n}}$$
$$b=-\frac{\bar k \rme^{\rmi \beta  \tau_0 } \left(1-\frac{\rmi \beta
   T}{n}\right)^{-n}}{2+\bar k \left(\frac{T \rme^{-\rmi \beta  \tau_0 }}{ 1+\frac{\rmi
   \beta  T}{n} }+\tau_0  \rme^{-\rmi \beta  \tau_0 }\right) \left(1+\frac{\rmi
   \beta  T}{n}\right)^{-n}}$$

We claim that  Re$a>0$ always holds. In fact, from (\ref{cegammas}), we have
$\rmi\beta+\rmi\omega_0+\Delta=\frac k {2}(1+\frac T
n\rmi\beta)^{-n}\rme^{-\rmi\beta\tau_0}$, then $$\begin{array}{ll}&\textrm{Sign~Re~}a\\&=\textrm{Sign~Re~}a^{-1}\\
&=\textrm{Sign~Re~}\frac{1+ \left(\frac{T  }{ 1+\frac{\rmi
   \beta  T}{n} }+\tau_0  \right) (\rmi\beta+\rmi\omega_0+\Delta)}{\rmi\beta+\rmi\omega_0+\Delta}\\
   &=\textrm{Sign~Re~}\left[\frac{1}{\rmi\beta+\rmi\omega_0+\Delta}+ \left(\frac{T  }{ 1+\frac{\rmi
   \beta  T}{n} }+\tau_0  \right) \right]\\&>0\end{array}$$
Thus only the two kinds of bifurcations shown in FIG. \ref{exa} can occur which are distinguished by the sign of Re$b$. Moreover, using the global Hopf bifurcation theorem \cite{Wuams,Wuams2}, we know all Hopf bifurcation branches are unbounded in the $k$ direction. When $k$ increases, after two, even more Hopf bifurcations occur, together with $\alpha=0$ is globally stable at $k=0$, we conclude the number of  coexisted coherent states gets larger.

   \subsection{The case without gap}
\begin{figure}[htbp]
  (a)~~~~~~~~~~~~~~~~~~~~~~~~~~~~~~~~~~~~~~~~~~~~~~~~~~~~~~~~~~~(b)\\
  \includegraphics[width=0.45\textwidth]{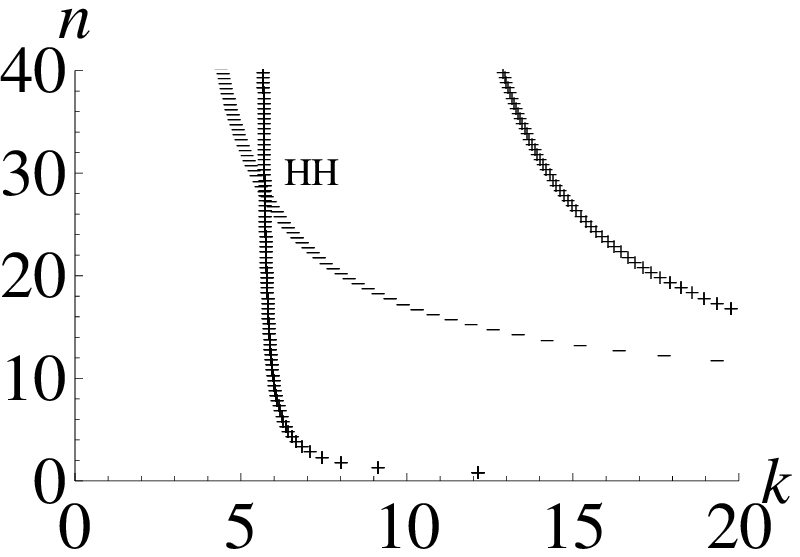}
 \includegraphics[width=0.45\textwidth]{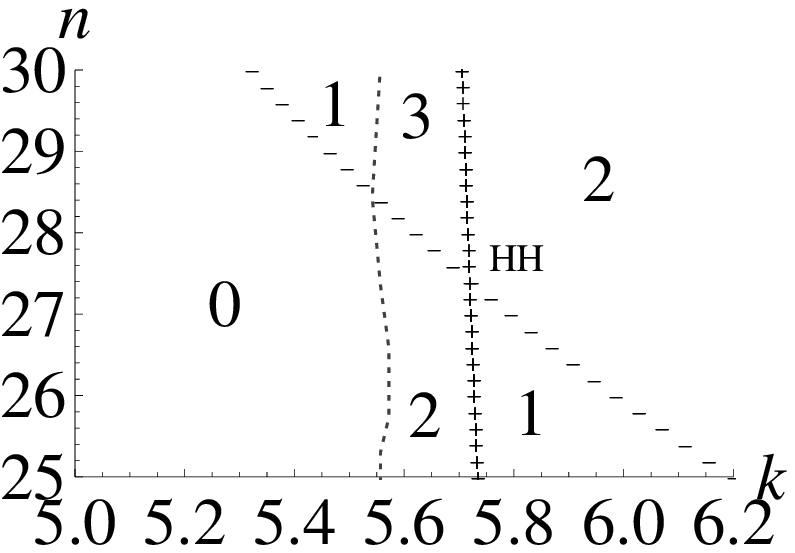}\\
  (c)~~~~~~~~~~~~~~~~~~~~~~~~~~~~~~~~~~~~~~~~~~~~~~~~~~~~~~~~~~~(d)\\
  \includegraphics[width=0.45\textwidth]{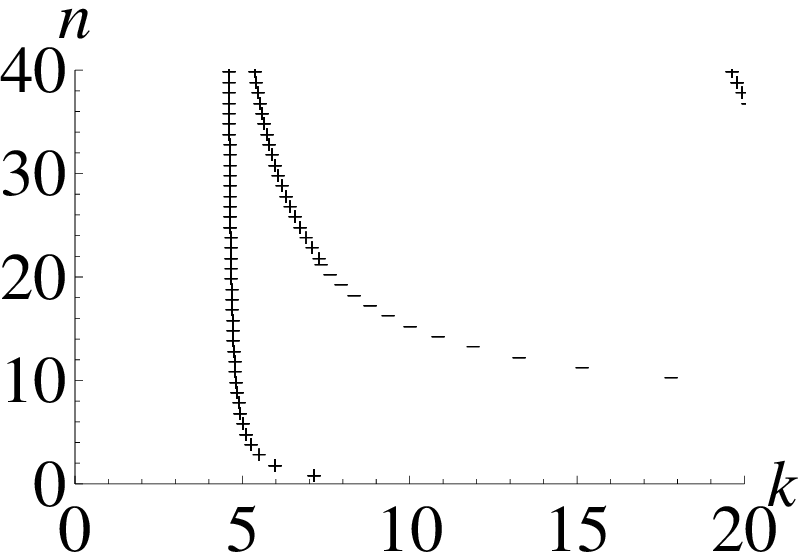}
 \includegraphics[width=0.45\textwidth]{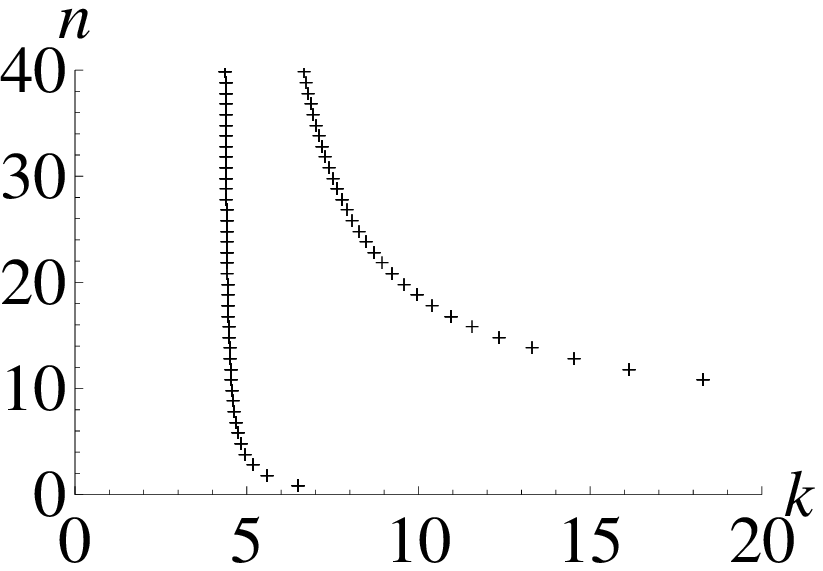}
  \caption{When $\omega_0=3$, $\Delta=1$, $T=3$, bifurcation curves are shown in (a), with a double Hopf bifurcation labeled by HH.  Near HH, number of coherent states is marked in (b). The dotted curve stands for the saddle-node bifurcation curve (sketched) illustrated in FIG. \ref{exa}. \lq\lq +/--" denotes subcritical/supercritical bifurcation. (c) and (d) are bifurcation curves for $T=1.15$ and $T=1$.}\label{tau0}
\end{figure}

   When $\tau_0=0$, i.e., $\tau$ degenerates into $\tilde\tau$ the Gamma distributed variable.  This is the case investigated by Lee \textit{et al} \cite{wai}. Hysteresis loop is observed when $T=3$, and the authors found smaller $n$ led to hysteresis loop while larger one did not. Using the method we established, we can calculate the bifurcation points and the signs of $b$ which are shown in FIG. \ref{tau0} (a) and (b). One   supercritical bifurcation curve  intersects with the other  subcritical bifurcation curve at a double Hopf bifurcation point HH. With the help of FIG. \ref{exa}, near the subcritical bifurcation we know there is a stable coherent state coexisting with the incoherent state, i.e., the hysteresis loop. When $T$ decreases, as shown in FIG. \ref{tau0} (c) and (d), supercritical bifurcation never occurs at the first bifurcation value $k$, which coincides with the previous results\cite{wai}. Finally, two remarks should be noticed that [i] in FIG. \ref{tau0} (b) the saddle node curve is a sketched one as we do not know how to calculate the exact values, theoretically; and [ii] near the double Hopf point HH, the dynamics may be more complicated such as the quasiperiodic behavior possibly existed  on 2-torus, even 3-torus\cite{Guckenheimer}.

   \subsection{The case with $\tau_0>0$ and  effect of  the gap}

\begin{figure}[htbp]
  (a)~~~~~~~~~~~~~~~~~~~~~~~~~~~~~~~~~~~~~~~~~~~~~~~~~~~~~~~~~~~(b)\\
  \includegraphics[width=0.45\textwidth]{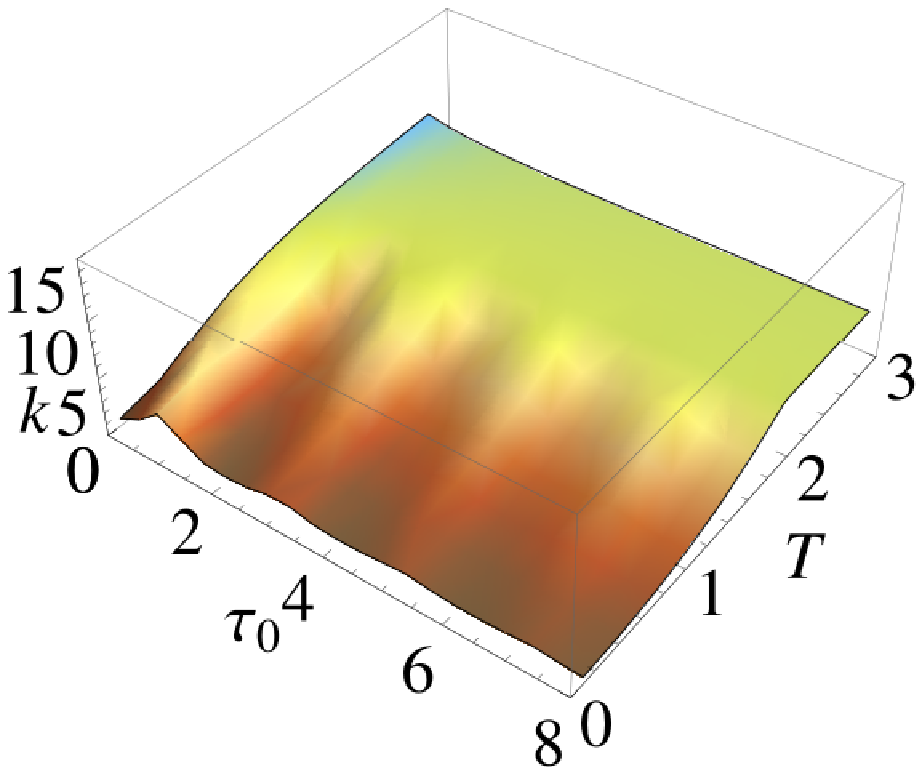}
 \includegraphics[width=0.45\textwidth]{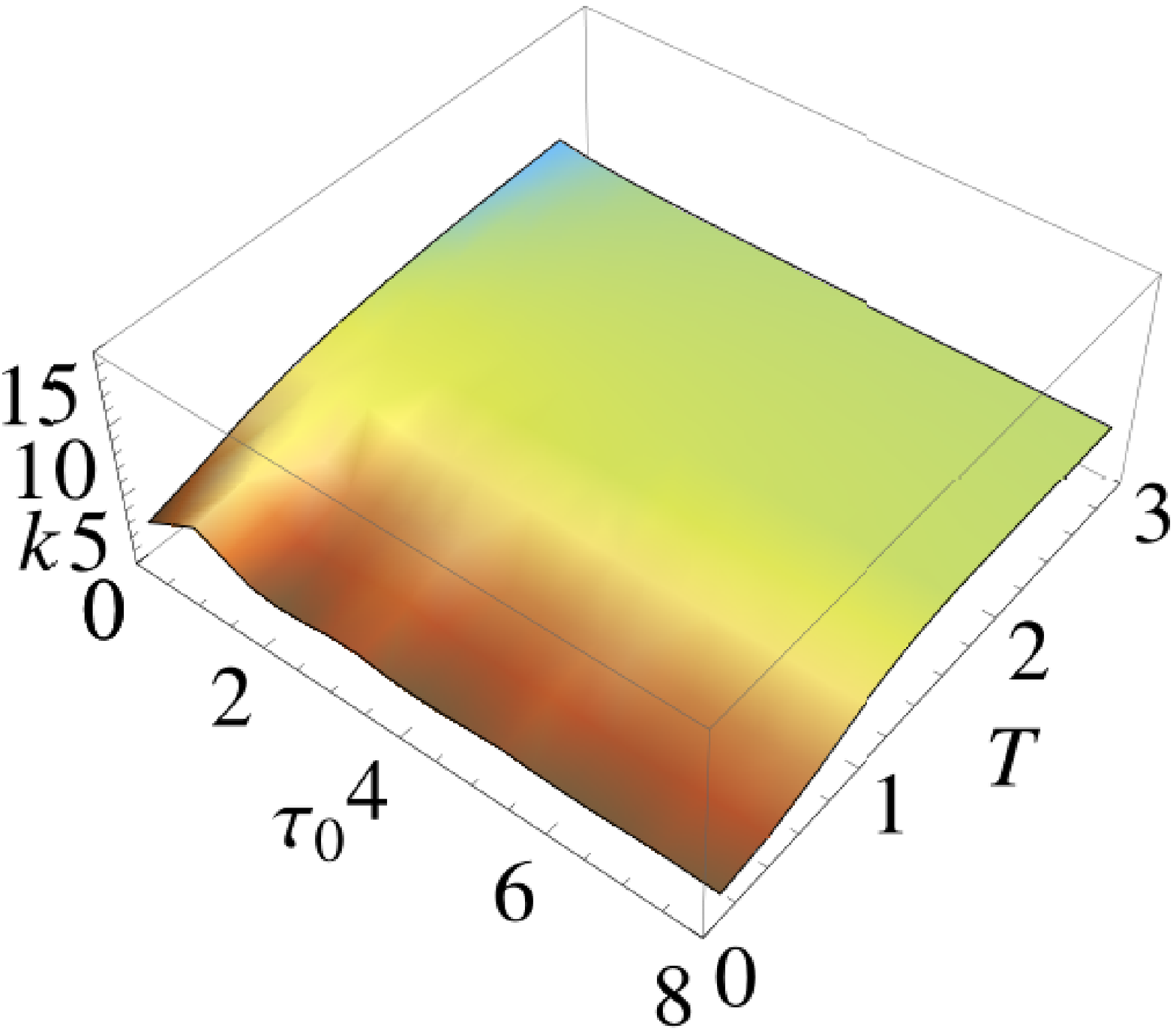}\\
   (c)~~~~~~~~~~~~~~~~~~~~~~~~~~~~~~~~~~~~~~~~~~~~~~~~~~~~~~~~~~~(d)\\
  \includegraphics[width=0.45\textwidth]{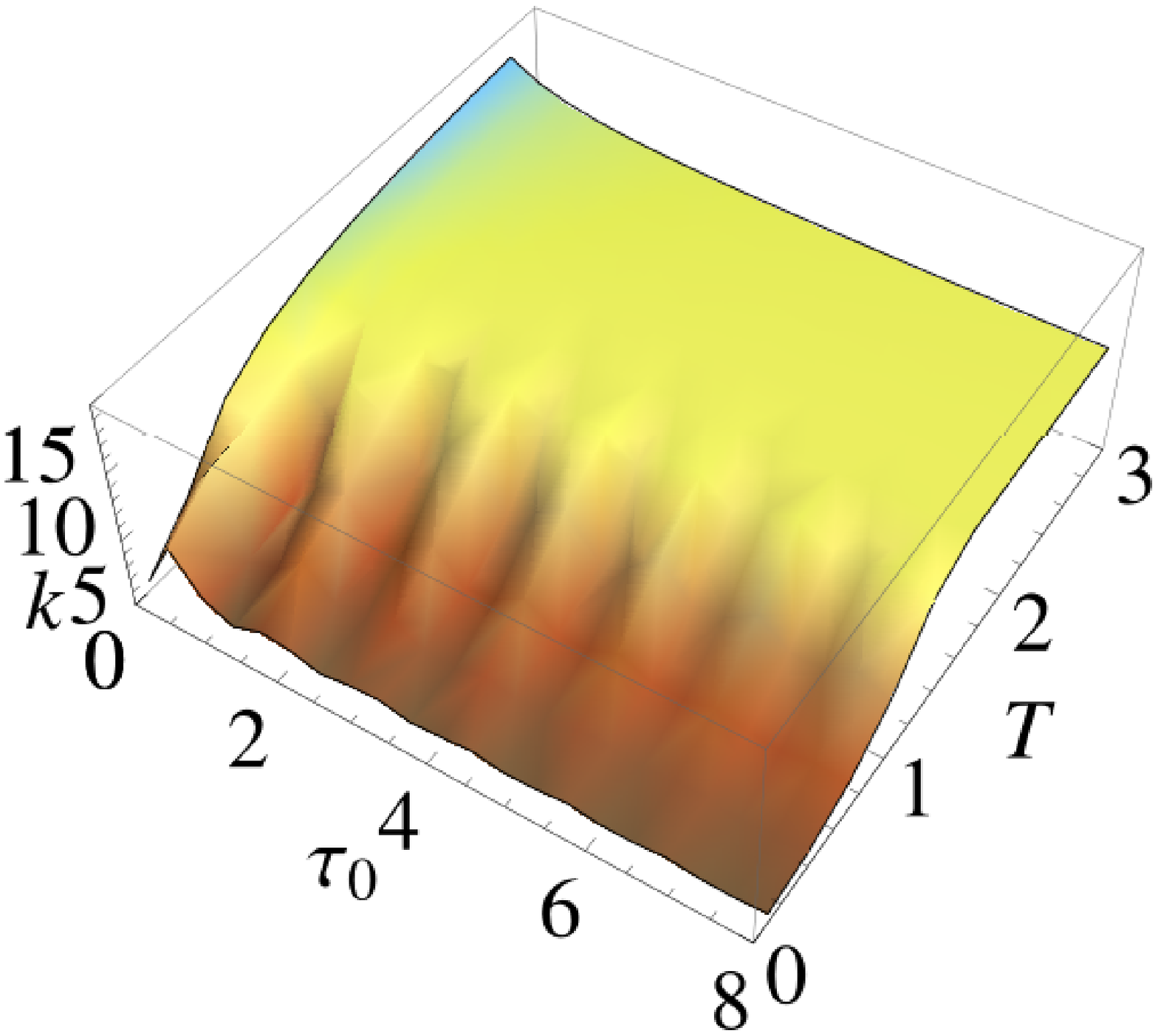}
 \includegraphics[width=0.45\textwidth]{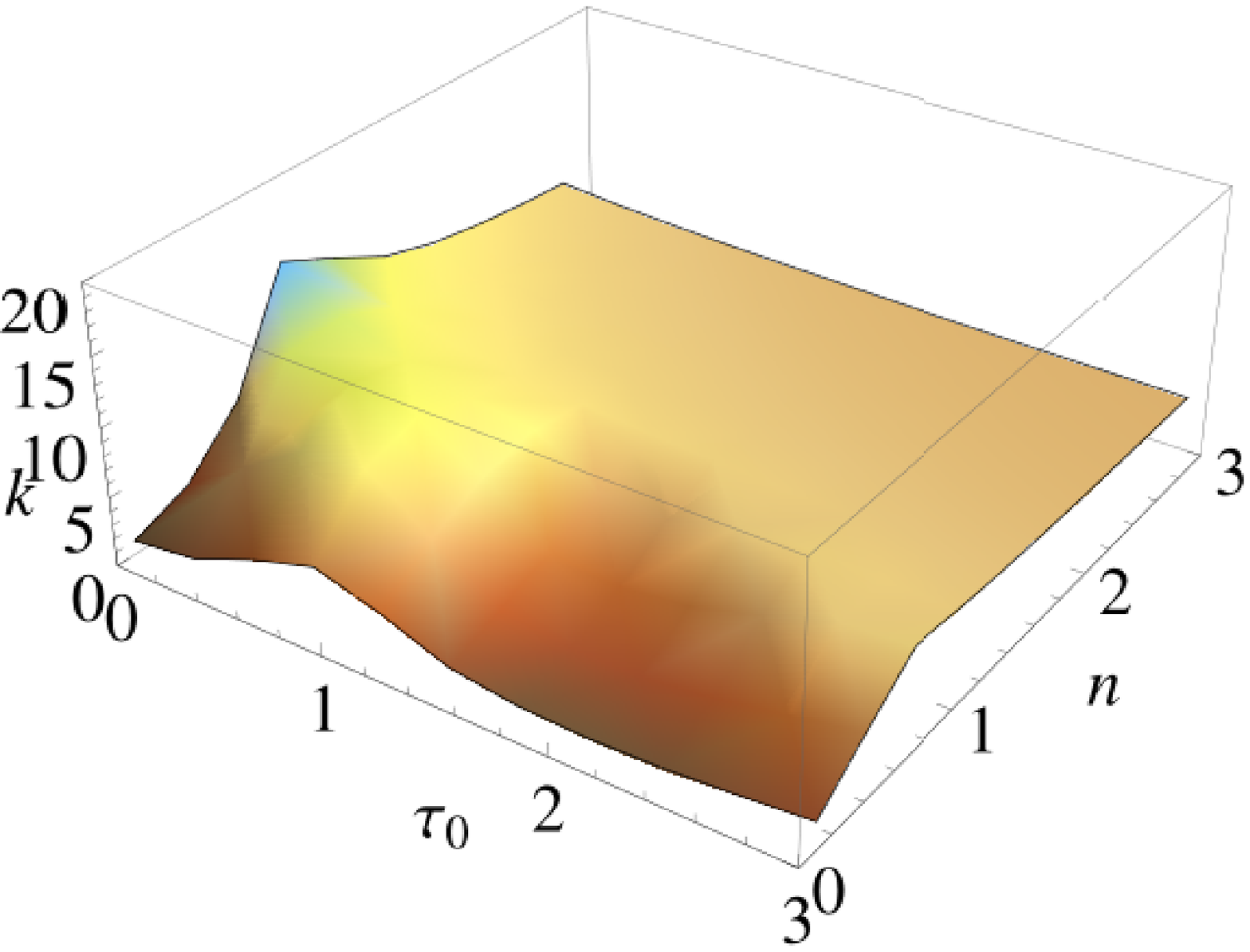}
  \caption{When   $n=3$, the first bifurcation value $k$ is calculated for $(T,\tau_0)\in[0,2.8]\times[0,8]$,    (a) $\omega_0=3$, $\Delta=0.3$;  (b) $\omega_0=3$, $\Delta=1$; (c) $\omega_0=5$, $\Delta=0.3$.  (d) When $\omega_0=3$, $\Delta=0.3$,  $T=3$, the first bifurcation value $k$ is calculated for $(n,\tau_0)\in[0,2.8]\times[0,3]$.}\label{tauT0}
\end{figure}

  In the rest part of this paper, we consider the case with a minimal positive delay, i.e., the case $\tau=\tilde\tau+\tau_0$ with $\tilde\tau$ a Gamma variable.   Using the method above, we can calculate the bifurcation values. In FIG. \ref{tauT0} (a)--(c), we find that increasing $\omega_0$ will delay the occurrence of bifurcation and increase in resonant structure\cite{ar9}
of the dependence of $k$ on $\tau_0$. When $T$, the mean  of Gamma distribution $\tilde\tau$ is larger, the effect of $\tau_0$ becomes weaker, thus the distributed delay $\tilde\tau$ acts predominantly. Increasing  the variance of the natural frequency $\Delta$ even further weakens the effect of the gap $\tau_0$.

\begin{figure}[htbp]
  (a)~~~~~~~~~~~~~~~~~~~~~~~~~~~~~~~~~~~~~~~~~~~~~~~~~~~~~~~~~~~(b)\\
  \includegraphics[width=0.45\textwidth]{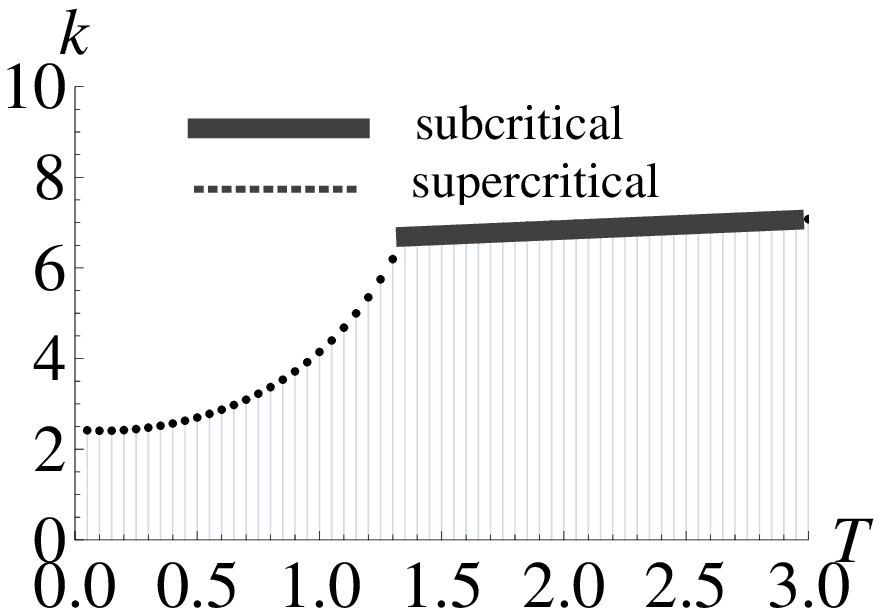}
 \includegraphics[width=0.45\textwidth]{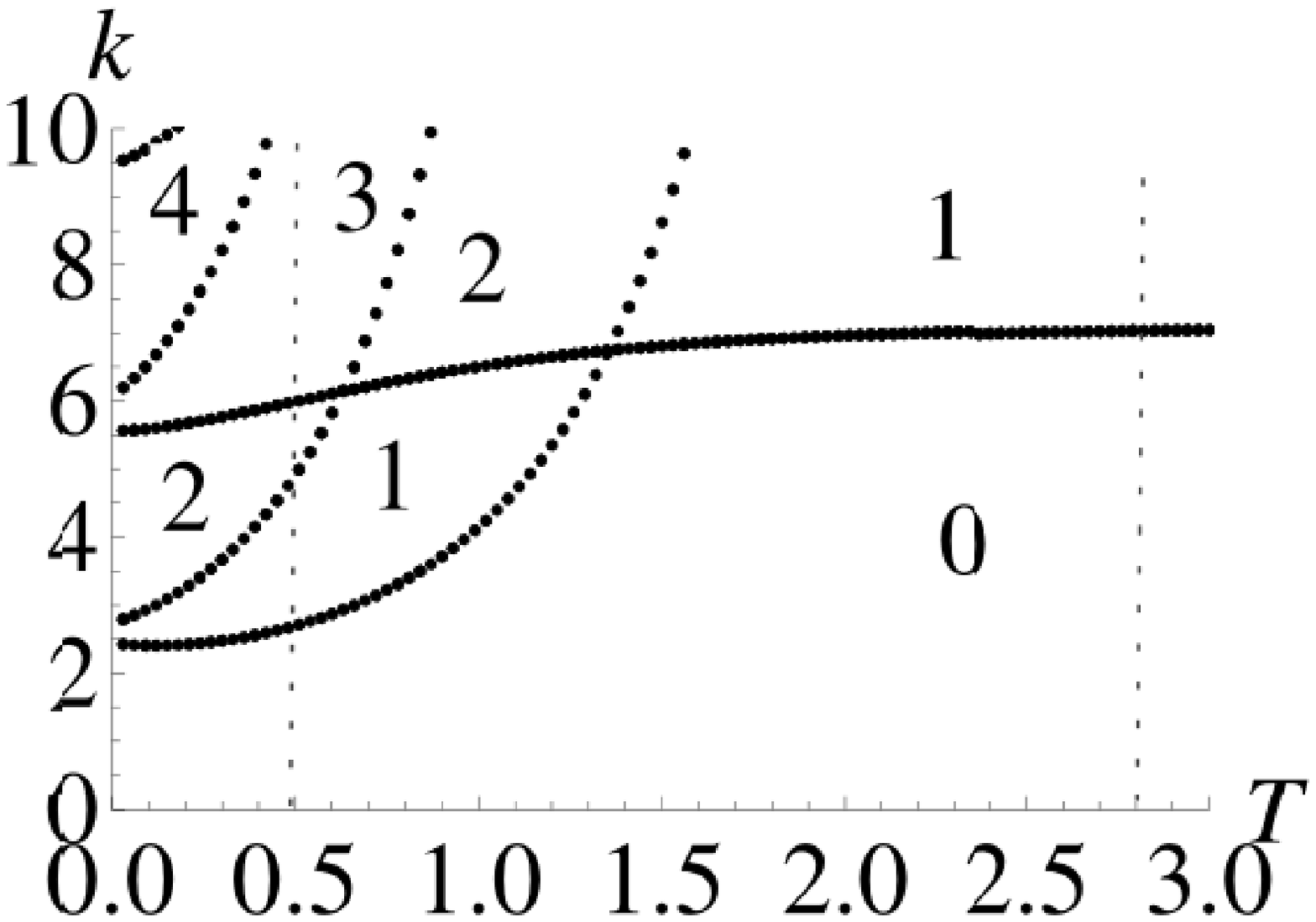} \\
  (c)~~~~~~~~~~~~~~~~~~~~~~~~~~~~~~~~~~~~~~~~~~~~~~~~~~~~~~~~~~~(d)\\
  \includegraphics[width=0.45\textwidth]{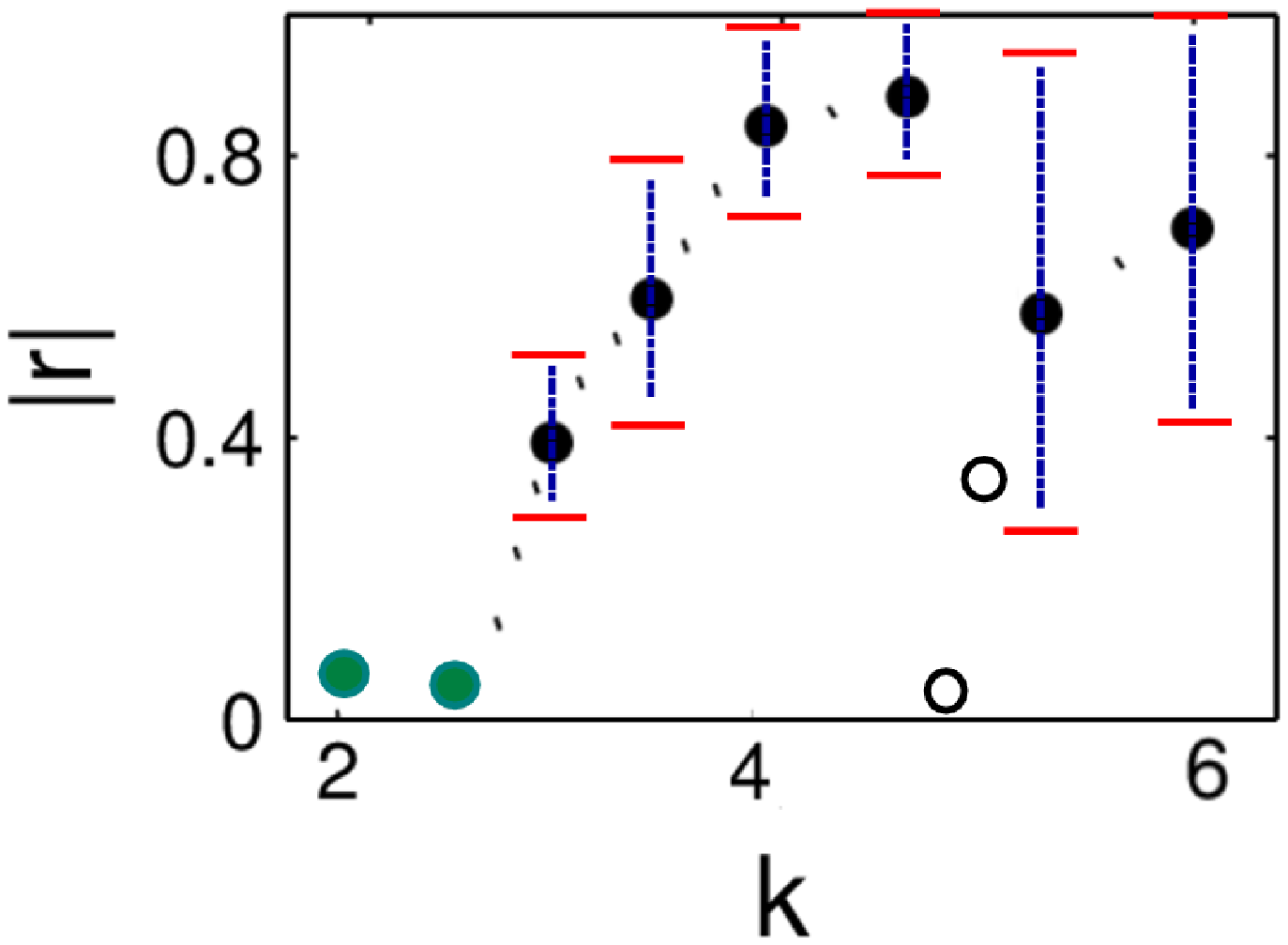}
 \includegraphics[width=0.45\textwidth]{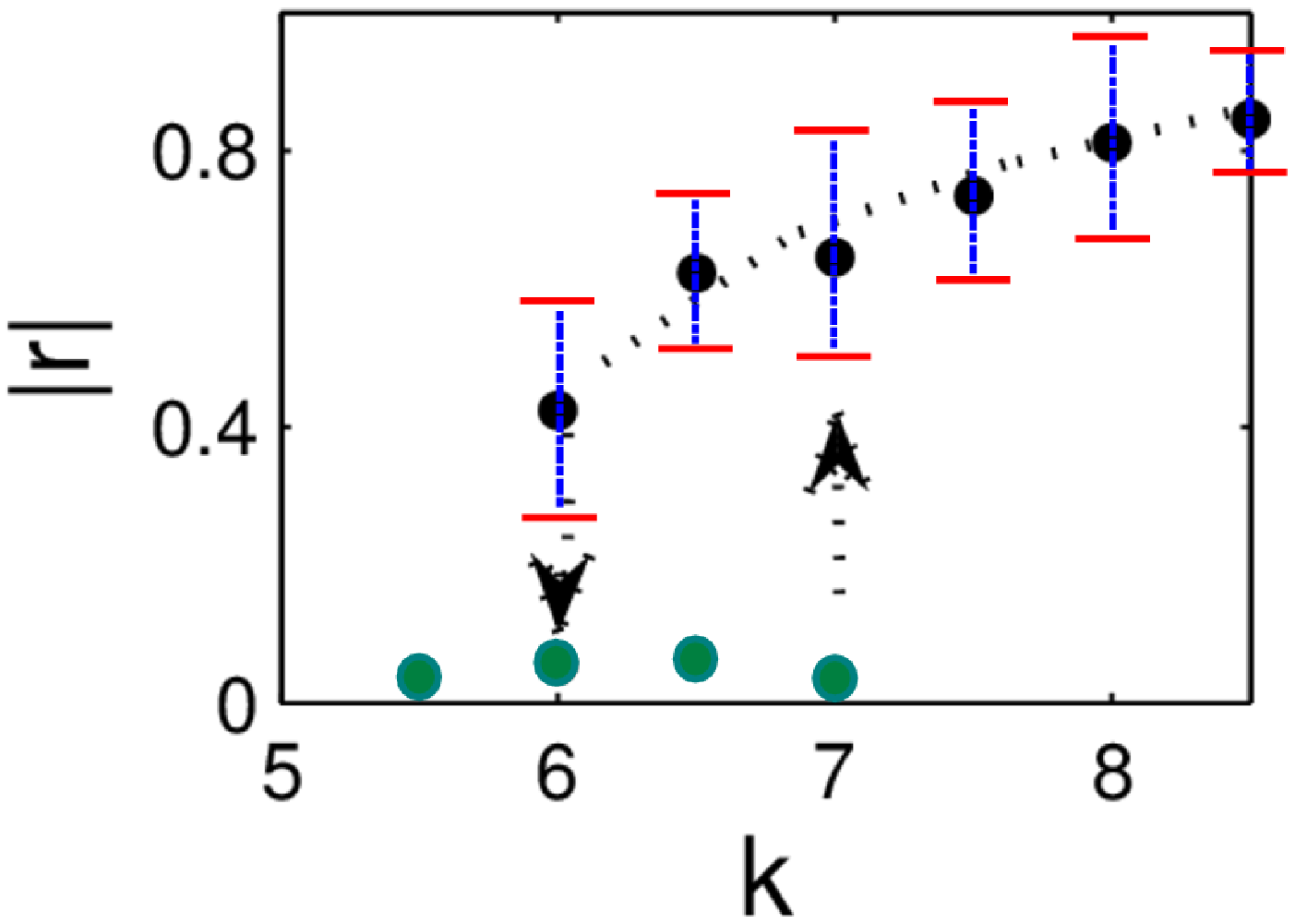}
  \caption{Fixing $\omega_0=3$, $\Delta=1$, $\langle\tilde\tau+\tau_0\rangle=T+\tau_0=3$, Var$(\tilde\tau+\tau_0)=\frac {T^2}n=3$,  the first bifurcation value and its direction is drawn in (a), and all bifurcation curves are illustrated in (b). The dark region stands for local stable region of incoherent state.  Numbers indicates the quantity of coherent attractors. When (c) $T=0.5$($\bar k=2.6992$) and (d) $T=2.8$($\bar k=7.0388$),   $|r|$ is shown in black dots by simulating system (\ref{model})  with $N=128$.}\label{w10}
\end{figure}
  By fixing $T$,   we are about to consider the interactional effect of  $\tau_0$ and the variance  of Gamma distribution characterized by $n$.
In FIG. \ref{tauT0} (d), we fix $T=3$ and investigate the effect of $\tau_0$ and $n$, where  we find increasing $n$   will decrease in resonant structure
of the dependence of $k$ on $\tau_0$, meanwhile weakens the effect of $\tau_0$.

Now we further discuss how the gap $\tau_0$ has an effect on the system dynamics,
when fixing the mean and variance of  $\tilde\tau+\tau_0$, i.e., the case shown in FIG. \ref{gammafig} (b).    Letting $\langle\tilde\tau+\tau_0\rangle=T+\tau_0=3$, Var$(\tilde\tau+\tau_0)=\frac {T^2}n=3$, the first Hopf bifurcation value and its direction are drawn in FIG. \ref{w10} (a). We find small  $T$ (i.e., $\langle\tilde\tau\rangle$) makes supercritical bifurcation occur at small $k$. When $T$ is large, the subcritical bifurcation occurs at large $k$, which means that in the case of fixed mean and variance of total delay, large proportion of the gap $\tau_0$ (i.e., small $T$) may destroy the hysteresis loop. In FIG. \ref{w10} (b), we draw all bifurcation curves and one can also find that larger  $\tau_0$   can significantly increase in the number of coherent attractors. When $T=0.5$ and $2.8$, respectively, simulations are carried out in FIG. \ref{w10} (c) and (d), where we find supercritical bifurcation and subcritical bifurcation in system (\ref{model}) with $128$ oscillators.

To give the simulations, we first use the method in \cite{ar9} to detect the stability of incoherent state by slightly perturbing the completely incoherent state: if the incoherent state is stable we show the minimal value of $|r|$ (green dots) among 100 times of simulations with random initial values; if the incoherent state is not stable we show the average value of $|r|$ (black dots) together with error bars (red segments) over 100 times of simulations with random initial values; otherwise if the incoherent state is not always stable, which is the case of coexistence, we show the average value of $|r|$ (black dots) together with error bars (red segments) over 100 times of simulations, too, but delete these values less than 0.2. Two black circles in FIG.  \ref{w10} (c) indicate speculative values of unstable coherent state bifurcating from $k=5$.
 
 As $k$ increases, two even more Hopf bifurcations  yield the coexistence of coherent states which makes the dynamical behavior more delicate.  As we have seen in FIG. \ref{w10} (c), when $k>5$, the standard variance of $|r|$ among 100 times of simulations is rather large, which reveals that several stable attractors coexist. They may be originated from the second branch of Hopf bifurcation (black circles). Another reason makes this standard variance large may be the unstable hysteresis loop near $k=6$, when a subcritical bifurcation occur. When $k$ is sufficiently near this critical point, hysteresis loop no longer contains stable states,    because after the first Hopf bifurcation the trivial solution of (\ref{sys2}) is certainly unstable. Notice that  this bifurcated unstable coherent state may turn into a stable one when $k$ increases far from the bifurcation point.  Nevertheless, the quantity of coherent states (stable or unstable) gets larger which makes the dynamics of Kuramoto model complicated.

The effect of the gap can also be explained in the following way: we notice that in case of Gamma distribution  larger $\tau_0$ means smaller $T$. As we have fixed the variance $\frac {T^2} n$, the excess kurtosis $\frac 6 n$ is larger. Notice that for Gamma distribution  kurtosis and skewness obey the same monotonic dependence, thus the kurtosis can be replaced by skewness. Since  excess kurtosis characterizes the  sharpness of the  peak of the distribution, in this point of view, larger kurtosis means that the sample data of delays are more ``concentrated'', which induces a supercritical bifurcation. Similarly we conclude that decentralized samples (smaller  excess kurtosis) of delay data induce subcritical bifurcations and  hysteresis loop. This result is summarized from calculation results (FIG. \ref{w10} (a) and (b)). However, giving an analytical relation between bifurcation properties and kurtosis needs more derivations. This is not an easy problem, thus is left as a further work.

\section{Conclusion}
In this paper, we establish a normal form method by extending  Nayfeh's  multiple scales to determine the properties of the bifurcated coherent states in a group of Kuramoto oscillators with heterogeneously  distributed   delays with a gap. Compared  with the previous work \cite{niuphyd}, where normal forms are derived by using a functional analysis method,  we find this a simple and useful way, on the Ott-Antonsen's manifold, to reveal the detailed dynamics near the critical values such as   the direction of bifurcation, stability of bifurcated coherent states and some coexistence phenomena. They can be determined by real part of two variables $a$ and $b$.

Some numerical results  indicate the  effect of the gap. As direct applications of our theory, how these parameters affect system dynamics is investigated.  For fixed   variance and expectation of total delay,  compared with the previous results\cite{wai}, we further find that larger gap $\tau_0$ (or larger excess kurtosis for Gamma distribution) has similar effect as the variance. It will (i) decrease the bifurcation values, (ii) induce a supercritical bifurcation hence avoid the hysteresis loop, and (iii) increase in the number of the coexisted coherent states.

\begin{acknowledgments}
The authors greatly appreciate the editor and the anonymous referees¡¯   comments and helpful suggestions which greatly improved the
presentation of the manuscript.
This research is supported by  National Natural Science Foundation of China under grant 11301117 and by Heilongjiang Provincial Natural Science Foundation under grant QC2014C003.
\end{acknowledgments}

\nocite{*}
\bibliography{aipsamp}% Produces the bibliography via BibTeX.

\begin{thebibliography}{99}
\bibitem{epsilon22} In fact, the traditional method of multiple scales requires $
\alpha(t;\epsilon):=\epsilon^{1/2}\alpha_1(T_0,T_1)+\epsilon \alpha_2(T_0,T_1)+\epsilon^{3/2}\alpha_3(T_0,T_1)+\cdots$. In this case, one can varify $\alpha_2=0$.
\bibitem{explain} When $k$ is in a sufficiently small neighbourhood of $\bar k$, the Hopf bifurcation value, the period of oscillation is given by $\frac {2\pi}{\beta}$. However, the amplitude of solutions decays very slow.
\bibitem {stable}  The stability can be obtained easily, by calculating the negative Floquet exponent at $|A(T_1)|=\sqrt{\frac{\nu\textrm{Re}a}{-\textrm{Re}b}}$.
\bibitem{Wuams2}In this case, with respect to a branch of Hopf bifurcating solutions, either it is  unbounded, or it contains $2m$ bifurcation points with Re$a>0$ at $m$ points and  Re$a<0$ at the rest $m$ points.
\end{thebibliography}

\end{document}